\input amstex
\documentstyle{amsppt}
\NoRunningHeads
\magnification=\magstep1
\baselineskip=12pt
\parskip=5pt
\parindent=18pt
\topskip=10pt
\leftskip=0pt
\pagewidth{32pc}
\pageheight{47pc}
\topmatter
\title Seshadri constant for a family of surfaces 
\endtitle
\author Keiji Oguiso  
\endauthor
\subjclass 
14C20, 14J10
\endsubjclass
\abstract 
The aim of this note is to study local and global Seshadri constants for a 
family of smooth surfaces with prescribed polarization. We shall first 
observe that given $\alpha$ being smaller than the square root of the degree 
of polarization, the set of local Seshadri constants in the range 
$(0, \alpha]$ is finite. This 
in particular implies that the square root of the degree of polarization is 
the only possible 
accumulation point of the set of local Seshadri constants. Next we shall 
remark the Zariski closedness of the set of points whose local Seshadri 
constants are in any given interval $(0, a]$. As applications, we shall 
also add a 
few remarks on the lower semi-continuity of both local and global Seshadri 
constants with respect to parameters involved, and on the minimality and the 
maximality of their infimum and supremum.
\endabstract

\endtopmatter

\document
\head 1. Introduction - Background and Results\endhead 

All the results in this note were entirely inspired by many interesting 
phenomena concerning Seshadri constants of algebraic surfaces being observed 
by Thomas Bauer ([B1], [B2], [B3]), and are nothing more than supplements for 
his works. New idea here is to study Seshadri constants for a 
family of surfaces, and perhaps, to regard Seshadri constants as a function 
of involved parameters a bit more consciously. 
\par
Throughout this note, we work over an algebraically closed field $k$ of 
any characteristic. A point means a closed point and a curve means an 
irreducible, reduced, complete curve, unless stated otherwise. 
A polarized surface is a pair of a smooth projective irreducible 
surface $X$ and an ample invertible sheaf $L$. 
By a family of polarized surfaces $(f : \Cal X \rightarrow \Cal B, \Cal L)$, 
we mean a proper flat morphism 
$f : \Cal X \rightarrow \Cal B$ over an irreducible (non-empty) noetherian 
scheme $\Cal B$ together with an $f$-ample invertible sheaf $\Cal L$ such that 
the fibers $(\Cal X_{t}, \Cal L_{t})$ ($t \in \Cal B$) are all polarized 
surfaces. Note that the intersection numbers $d := (\Cal L_{t}^{2})$ are 
independent of $t \in \Cal B$. We call this integer the degree of 
$(f : \Cal X \rightarrow \Cal B, \Cal L)$. Note also that a polarized 
surface is nothing but ``a family of polarized surfaces over 
$\text{Spec}\, k$''. 
\par
Let $(X, L)$ be a polarized surface of degree $d$. As well-known, for a given 
point $x \in X$, the local Seshadri constant 
$\epsilon(L, x)$ of $(X, L)$ at $x$ is defined to be $\epsilon(L, x) := 
\text{inf}_{x \in C} (L.C)/m_{x}(C)$, where $m_{x}(C)$ is the multiplicity of 
$C$ at $x$ and the infimum is taken over all curves $C \subset X$ 
passing through $x$. It is also well-known that $\epsilon(L, x) = 
\text{max}\{s \in \Bbb R \vert \pi^{*}L - sE \, 
\text{is nef.} \}$, where $\pi : \tilde{X} \rightarrow X$ is the blow up at 
$x \in X$ and $E$ is the exceptional curve. Therefore, as observed 
by A. Steffens [S, Proposition 4], the real version of the Nakai-Moishezon 
criterion [CP] implies $0 < \epsilon(L, x) \leq \sqrt{d}$, and 
$\epsilon(L, x) \in \Bbb Q$ unless $\epsilon(L, x) = \sqrt{d}$ and 
$\sqrt{d} \not\in \Bbb Z$. (See also Lemma (2.1).) The global Seshadri 
constant $\epsilon(L)$ is defined to 
be $\epsilon(L) := \text{inf}_{x \in X} \epsilon(L, x)$ 
([EL], [B3, Section 1]). Then, one has also $0 < \epsilon(L) \leq \sqrt{d}$, 
where the first inequality is because of the well-known criterion of ampleness 
due to Seshadri. 
\par

The first aim of this note is to observe the following finiteness:

\proclaim{Theorem 1} Let $(f : \Cal X \rightarrow \Cal B, \Cal L)$ be a 
family of polarized surfaces of degree $d$ and set $\Sigma := 
\{\epsilon(\Cal L_{t}, x_{t}) \vert t \in \Cal B, x_{t} \in X_{t}\}$.
Then, for each given number $\alpha \in \Bbb R$ such 
that $\alpha < \sqrt{d}$, the set $\Sigma \cap (0, \alpha]$ is finite. 
In particular, the only 
possible accumulation point of the set $\Sigma$ is $\sqrt{d}$ (and the 
accumulation would be, of course, from the below if it happens). 
\endproclaim 

This result was inspired by [B3, Theorem 2.1], in which very precise 
possiblities of values $\epsilon(L, x)$ belonging to the range $(0, 2)$ are 
given for surfaces $(X, \Cal O_{X}(1))$ embedded into projective spaces by 
$\Cal O_{X}(1)$. 
\par

Combining this with Lemma (2.1), one immediately obtains: 

\proclaim{Corollary 2} Let $(X, L)$ be a polarized surface of degree 
$d$. Then, $\epsilon(L) \in \Bbb Q$ unless $\epsilon(L) = 
\sqrt{d}$ and $\sqrt{d} \not\in \Bbb Z$. In particular, $\epsilon(L) \in 
\Bbb Q$ if $\sqrt{d} \in \Bbb Z$. More precisely, there exists a point 
$x \in X$ and a curve $x \in C \subset X$ such that 
$\epsilon(L) = \epsilon(L, x) = (L.C)/m_{x}(C)$ unless 
$\epsilon(L) = \sqrt{d}$. In particular, 
there always exists a point $x \in X$ such that $\epsilon(L) = 
\epsilon(L, x)$, i.e. one always has $\inf_{x \in X} \epsilon(L, x) 
= \min_{x \in X} \epsilon(L, x)$. \qed 
\endproclaim 

Note that in [B1] and [B2, Appendix by T. Bauer and T. Szemberg], the 
rationality of $\epsilon(L)$ for quartic K3 surfaces and for polarized 
abelian surfaces is proved.
\par

The second aim of this note is to observe the following closedness:

\proclaim{Theorem 3} Let $(f : \Cal X \rightarrow \Cal B, \Cal L)$ be 
a family of polarized surfaces of degree $d$ 
and set $\Cal X(a) := \{x_{t} \vert 
t \in \Cal B, x_{t} \in \Cal X_{t}, \epsilon(\Cal L_{t}, x_{t}) \leq a\} 
(\subset \Cal X)$, where $a$ is any real number. Then, $\Cal X(a)$ is Zariski 
closed in $\Cal X$. 
\endproclaim 

Theorem 3 together with an observation due to L. Ein and R. Lazarsfeld 
[EL, Theorem], in particular, implies the following slight refinement of 
their result: 

\proclaim{Corollary 4} 
Given a positive number $0 <\delta \in \Bbb R$, the set $\{x \in X \vert 
\epsilon(L, x) \leq 1 -\delta \}$ is finite for each polarized surface 
$(X, L)$. \qed  
\endproclaim 

The following semi-continuity is now also clear by Theorems 1 and 3: 

\proclaim{Corollary 5}  
\roster
\item 
For each fixed $t \in \Cal B$, the function $y = \epsilon(x) :=
\epsilon (\Cal L_{t}, x)$ of $x \in \Cal X_{t}$ 
is lower semi-continuous with respect to the $I$-topology of $\Cal X_{t}$. 
In addition, if $\epsilon(x) < \sqrt{d}$ at a point $x \in \Cal X_{t}$, 
then $y = \epsilon(x)$ is lower semi-continuous at this $x \in \Cal X_{t}$ 
also in the Zariski topology of $\Cal X_{t}$. 
\item 
The function $y = \epsilon(t) :=
\epsilon (\Cal L_{t})$ of $t \in \Cal B$ , where 
$\epsilon (\Cal L_{t})$ is the 
global Seshadri constant of $(\Cal X_{t}, \Cal L_{t})$, is lower 
semi-continuous with respect to the $I$-topology of $\Cal B$. In addition, 
if $\epsilon(t) < \sqrt{d}$ at a point $t \in \Cal B$, 
then $y = \epsilon(t)$ is lower semi-continuous at this $t \in \Cal B$ 
also in the Zariski topology of $\Cal B$. 
\endroster
Here, by the $I$-topology of a noetherian scheme $S$, we mean a topology of 
$S$ in which the open sets are $\emptyset$, $S$, and $S - T$, where $T$ is 
a union of at most countably many closed subschemes of $S$, and a real 
valued function $y = F(x)$ on a topological space $S$ is said to be 
lower semi-continuous at a point $x \in S$ if there exists an open 
subset $U (\subset S)$ such that $x \in U$ and that $F(x) \leq F(x')$ for all 
$x' \in U$. 
\qed 
\endproclaim 

This result was much inspired by work [B1] on the global Seshadri constants of 
quartic surfaces: They are mostly constant but jump below at special locus in 
the moduli. 
\par
In the light of Theorem 1, the values $\sigma(\Cal L_{t}) := 
\sup_{x_{t} \in \Cal X_{t}}\{\epsilon(\Cal L_{t}, 
x_{t})\}$ for each $t \in \Cal B$ and $\sigma(\Cal L) := 
\sup_{t \in \Cal B} \{\sigma(\Cal L_{t})\}$ 
might be of some interest. Concerning these values, combining our Theorems 
together with the fact that any union of countably many proper closed subsets 
does not cover the whole irreducible scheme if the base field is uncountable, 
one can immediately obtain the following: 

\proclaim{Corollary 6} Assume that the base field is uncountable. 
Then: 
\roster 
\item 
There exists a dense subset $U_{t} \subset \Cal X_{t}$ 
for each $t \in \Cal B$ such that $\sigma(\Cal L_{t}) = 
\epsilon(\Cal L_{t}, x_{t})$ for all $x_{t} \in U_{t}$;
\item 
There exists a dense subset $\Cal U \subset \Cal B$ such that 
$\sigma(\Cal L) = \sigma(\Cal L_{t})$ for all $t \in \Cal U$.
\endroster
In particular, there also exist 
$t \in \Cal B$ and $x_{t} \in \Cal X_{t}$ such that 
$\sigma(\Cal L) = \epsilon(\Cal L_{t}, x_{t})$, 
again both supremum is maximum,  
and $\sigma(\Cal L)$ and $\sigma(\Cal L_{t})$ are all rational unless 
$\sigma(\Cal L) = \sqrt{d}$ and $\sqrt{d} \not\in \Bbb Z$. \qed
\endproclaim 

There seems to be no known examples of $(X, L)$ of degree $d$ such that 
$\sqrt{d} \not\in \Bbb Z$ but $\epsilon(L, x) = \sqrt{d}$ at some $x \in X$. 
\par

All of the statements are standard applications of the existence 
of the relative Hilbertscheme (eg. [K]) together with the well-known 
Lemma (2.1) below. 

\head Acknowledgement \endhead

This note has been written up during the author's stay in Universit\"at Essen 
under the financial support by the Alexander-von-Humboldt fundation. First of 
all, the author would like to express his sincere thanks to Professors 
H\'el\`ene Esnault and Eckart Viehweg and the Alexander-von-Humboldt fundation 
for making the author's stay possible. The author would like to express his 
best thanks to Tomasz Szemberg for valuable discussions on the Seshadri 
constant as well as several helpful comments on this note, and again to 
Professor Eckart Viehweg for his suggestion to study 
the Seshadri constants for a family, without either of which the author could 
not carry out this work. Last but not least at all, the author would like to 
express his hearty thanks to Professor Yujiro Kawamata and Miss Kaori Suzuki 
for their warm encouragements through e-mails.  

\head 2. Proof of Theorems \endhead 

The following easy but remarkable Lemma is well-known ([CP], 
[S, Proposition 4]): 

\proclaim{Lemma (2.1)} Let $(X, L)$ be a polarized surface of degree $d$. 
If $\epsilon(L, x) < \sqrt{d}$, 
then there exists  a curve $x \in C \subset X$ such that 
$\epsilon(L, x) = (L.C)/m_{x}(C)$. In particular, 
$\epsilon(L, x) \in \Bbb Q$ unless $\epsilon(L, x) = \sqrt{d}$ and 
$\sqrt{d} \not\in \Bbb Z$.\qed
\endproclaim 

\proclaim {Lemma (2.2)} Let $(f : \Cal X \rightarrow \Cal B, \Cal L)$ be a 
family of polarized surfaces of degree $d$. Let $a \in \Bbb Q$ such that 
$0 < a < \sqrt{d}$. 
Then, there exists an integer $B := B(a)$, depending only on $a$, such that 
$(\Cal L_{t}.C_{t}) \leq B$ for any points $t \in \Cal B$, $x_{t} 
\in \Cal X_{t}$ 
and for any curve $x_{t} \in C_{t} \subset \Cal X_{t}$ with 
$(\Cal L_{t}.C_{t})/m_{x_{t}}(C_{t}) \leq a$. 
\endproclaim 

\remark{Remark} The idea of proof below was much inspired by [B3] and is 
indeed nothing more than a simple modification of arguments there towards our 
aim. 
\endremark

\demo{Proof} Since the statement for $\Cal L$ follows from the one for 
$\Cal L^{\otimes l}$ for a positive integer $l$, by the Serre vanishing 
Theorem, we may assume that $R^{i}f_{*}\Cal L^{\otimes n} 
= 0$ for all $i > 0$ and $n > 0$. Then, $h^{i}(\Cal L_{t}^{\otimes n}) 
= 0$ for all $t \in \Cal B$ as well. Therefore, by the Riemann-Roch formula, 
one has 
$h^{0}(\Cal L_{t}^{\otimes n}) = \chi(\Cal L_{t}^{\otimes n}) 
= n^{2}d/2 + nc/2 + c'$, where $c$ and $c'$ are integers independent 
of $t \in \Cal B$.   
Since $a \in \Bbb Q$, there exists a sequence of integers such that 
$n_{k} > 0$, $n_{k}a \in \Bbb Z$ and that $\text{lim}_{k \rightarrow 
\infty} n_{k} = \infty$. Set $l(n_{k}) := h^{0}(\Cal L_{t}^{\otimes n_{k}}) - 
(n_{k}a+2)(n_{k}a+1)/2$. 
Then, $l(n_{k}) = (d - a^{2})n_{k}^{2}/2 + (c-3a)n_{k}/2 + (c'-1)$. 
Since $d - a^{2} > 0$, one has 
$l(n_{k}) > 0$ for $k$ being large. Set $M := n_{k}$ for one of 
such $n_{k}$. Then, for any $x_{t} \in \Cal X_{t}$, there exists an effective 
divisor $D_{t}$ on $\Cal X_{t}$ such that  
$D_{t} \in \vert \Cal L_{t}^{\otimes M} \vert$ and that 
$m_{x_{t}}(D_{t}) \geq Ma + 1$ by 
$l(M) > 0$. Let $C_{t} \subset \Cal X_{t}$ be a curve such that 
$x_{t} \in C_{t}$ and that 
$(\Cal L_{t}.C_{t})/m_{x_{t}}(C_{t}) \leq a$. If 
$\text{Supp}(C_{t}) \not\subset \text{Supp}(D_{t})$, then, by the 
irreducibility 
of $C_{t}$, one would obtain $(\Cal L_{t}^{\otimes M}.C_{t}) = 
(D_{t}.C_{t}) \geq 
m_{x_{t}}(D_{t})\cdot m_{x_{t}}(C_{t}) \geq (Ma+1)
(\Cal L_{t}.C_{t})/a 
> (\Cal L_{t}^{\otimes M}.C_{t})$, a contradiction. 
Thus, $\text{Supp}(C_{t}) 
\subset \text{Supp}(D_{t})$. Since $C_{t}$ is also reduced and since 
$\Cal L_{t}$ 
is ample on $\Cal X_{t}$, one then obtains $(\Cal L_{t}.C_{t}) \leq 
(\Cal L_{t}.D_{t}) = Md$. 
Therefore, $B := Md$ provides a desired integer. \qed 
\enddemo

\proclaim{Lemma (2.3)} Let $a_{i}$, $b_{i}$ ($1 \leq i \leq n$) be positive 
real numbers. Then, $\min_{i}\{a_{i}/b_{i}\} \leq 
\sum_{i=1}^{n} a_{i}/\sum_{i=1}^{n} b_{i} \leq \max_{i}\{a_{i}/b_{i}\}$. 
\endproclaim 

\demo{Proof} Induction on $n$ plus elementary calculation. \qed 
\enddemo

\head Proof of Theorem 1 \endhead  
Since $\epsilon(\Cal L_{t}, x_{t}) = 
\epsilon((\Cal L^{\otimes n})_{t}, x_{t})/n$ for any $t \in \Cal B$, 
$x_{t} \in \Cal X_{t}$ and for any positive integer $n$, 
we may assume that $\Cal L$ is $f$-very ample. 
Set $\Cal S := \{\epsilon(\Cal L_{t}, x_{t}) \vert 
t \in \Cal B, x_{t} \in \Cal X_{t}\} \cap (0, \alpha]$. If 
$\Cal S = \emptyset$, then the result is true. Therefore we may assume that 
$\Cal S \not= \emptyset$. 
Let $s = \epsilon(\Cal L_{t}, x_{t}) \in \Cal S$. Since $\alpha < 
\sqrt{d}$, by Lemma (2.1), there exists a curve $C_{t} \subset \Cal X_{t}$ 
such that $x_{t} \in C_{t}$ 
and that $s = (\Cal L_{t}.C_{t})/m_{x_{t}}(C_{t})$. Since there is a 
rational number $a$ 
such that $\alpha < a < \sqrt{d}$, by Lemma (2.2), there exists an 
integer $B$ (independent of $s \in \Cal S$) such that 
$(\Cal L_{t}.C_{t}) \leq B$ for all such pairs 
$(x_{t} \in C_{t})$ above. Since 
$\Cal L_{t}$ is very ample on $\Cal X_{t}$, for each such $x_{t} \in 
C_{t}$, there exists an element $D_{t} \in \vert \Cal L_{t} \vert$ 
such that 
$x_{t} \in D_{t}$ but $y \not\in D_{t}$ for some $y \in C_{t}$. 
Since $C_{t}$ 
is irreducible, $C_{t}$ and $D_{t}$ then meet properly. Therefore, 
by $x_{t} \in D_{t}$, we calculate 
$1 \leq m_{x_{t}}(C_{t}) \leq m_{x_{t}}(D_{t})\cdot m_{x_{t}}
(C_{t}) \leq 
(D_{t}.C_{t}) = (\Cal L_{t}.C_{t}) \leq B$. Since 
$m_{x_{t}}(C_{t})$ 
and $(\Cal L_{t}.C_{t})$ are integers, the possible pairs of values
$(m_{x_{t}}(C_{t}), (\Cal L_{t}.C_{t}))$ are then finite. 
Therefore, $\Cal S$ is finite as well. \qed

\head Proof of Theorem 3 \endhead 

By Theorem 1, we may assume that $0 < a < \sqrt{d}$ and $a \in \Bbb Q$. 
Let $x_{t} \in \Cal X(a)$, where we denote $t := f(x_{t}) \in \Cal B$. 
Since $a < \sqrt{d}$, by Lemma (2.1), there exists a curve 
$x_{t} \in C_{t} \subset \Cal X_{t}$ such that $\epsilon(\Cal L_{t}, x_{t}) = 
(\Cal L_{t}.C_{t})/m_{x_{t}}(C_{t})$. Let us consider the 
product $\Cal H^{0} 
\times_{\Cal B} \Cal H^{1}$ of the relative Hilbertschemes of 
points $\Cal H^{0}$ and the relative Hilbertschemes of 
one dimensional subschemes $\Cal H^{1}$ of our family 
$(f : \Cal X \rightarrow \Cal B, \Cal L)$ 
and denote by $\Cal K(a)$ the subset consisting of all 
$[x_{t} \in C_{t}]$ as above (here, $t \in \Cal B$ also varies). 
Define $\Cal H(a)$ to be the Zariski closure of $\Cal K(a)$ in $\Cal H^{0} 
\times_{\Cal B} \Cal H^{1}$. By Lemma (2.2), there exists 
an integer $B$ such that $(\Cal L_{t}.C_{t}) \leq B$ for all 
$[x_{t} \in C_{t}] \in 
\Cal K(a)$. Therefore, the number of the irreducible components of the 
relative Hilbertscheme meeting $\Cal K(a)$ are then finite. Thus, 
$\Cal H(a)$ has also finitely many irreducible components. 
Let $\Cal H(a)_{i}$ ($1 \leq i \leq I$) be all the irreducible components of 
$\Cal H(a)$ and $\Cal C(a)_{i} \rightarrow \Cal H(a)_{i}$ be the universal 
family. Note that the natural morphism 
$\Cal C(a)_{i} \rightarrow \Cal H(a)_{i} \rightarrow \Cal B$ is projective. 
Write $\Cal K(a)_{i} = \Cal K(a) \cap \Cal H(a)_{i}$. Then 
$\Cal H(a)_{i}$ is the Zariski closure of $\Cal K(a)_{i}$. Take 
$[x_{0} \in C_{0}] \in \Cal H(a)_{i}$ and put $t_{0} := f(x_{0}) 
(= f(C_{0}))$. Note that here $C_{0}$ might be neither irreducible nor 
reduced, but is certainly Cartier on $\Cal X_{t_{0}}$ and then has no 
embedded points, because of the universal closedness of the relative 
Cartier divisor functor for $f : \Cal X \rightarrow \Cal B$ being smooth 
([K, Page 18 Theorem 1.13]). 
Let $r(x, y, h) = 0$ 
be the local equations of the pointed curves $P(h) \in C(h)$ in $\Cal X$ 
such that $[P(h) \in C(h)] \in \Cal U (\subset \Cal H(a)_{i})$, where 
$\Cal U$ is a neighbourhood of $[x_{0} \in C_{0}]$, $(x, y)$ are 
fiber coordinates of $f$ around $x_{0}$ and $h$ stands for the 
parameters of $\Cal U (\subset \Cal H(a)_{i})$. Write $P(h) = (x(h), y(h))$. 
Then, for any given $m$, the locus such that $m_{P(h)}(C(h)) \geq m$ is 
defined by 
the vanishing of all the coefficients of terms of order $\leq (m-1)$ with 
respect to $x - x(h)$, $y - y(h)$ of 
$r(x, y, h)$. Thus, $\{[y \in D] \in \Cal H(a)_{i} \vert 
m_{y}(D) \geq m \}$ is Zariski closed in $\Cal H(a)_{i}$. Set 
$M := \min \{m_{x}(C) \vert [x \in C] \in \Cal K(a)_{i}\}$ and take 
$[x' \in C'] \in \Cal K(a)_{i}$ such that $M = m_{x'}(C')$. 
We set $t' := f(x')$.  
Then, the set $\Cal N(a)_{i} := \{[y \in D] \in \Cal H(a)_{i} \vert m_{y}(D) 
\geq M \}$ is Zariski closed in $\Cal H(a)_{i}$ and contains 
$\Cal K(a)_{i}$. Since $\Cal H(a)_{i}$ was the Zariski closure of 
$\Cal K(a)_{i}$, we have then $\Cal H(a)_{i} = \Cal N(a)_{i}$. In particular, 
$m_{x_{0}}(C_{0}) \geq M$. Since$(\Cal L_{t_{0}}.C_{0}) = (\Cal L_{t'}.C')$, 
we obtain $(\Cal L_{t_{0}}.C_{0})/m_{x_{0}}(C_{0}) \leq 
(\Cal L_{t_{0}}.C_{0})/M = (\Cal L_{t'}.C')/m_{x'}(C')$. Combining this with 
$(\Cal L_{t'}.C')/m_{x'}(C') \leq a$ (the definition of $\Cal K(a)_{i}$), 
one obtains 
$(\Cal L_{t_{0}}.C_{0})/m_{x_{0}}(C_{0}) \leq a$ as well.  
Let $C_{0} = 
\sum_{j} a_{j}E_{j} + \sum_{l} b_{l}F_{l}$ be the irreducible 
decomposition of $C_{0}$ such that $x_{0} \in E_{j}$ but 
$x_{0} \not\in F_{l}$. Since $\Cal L_{t_{0}}$ is ample, one has 
$(\Cal L_{t_{0}}.C_{0}) = \sum_{j} a_{j}(\Cal L_{t_{0}}.E_{j}) + 
\sum_{l} b_{l}(\Cal L_{t_{0}}.F_{l}) 
\geq \sum_{j} a_{j}(\Cal L_{t_{0}}.E_{j})$. One also has
$m_{x_{0}}(C_{0}) = \sum_{j}a_{j}m_{x_{0}}(E_{j})$. 
Then, by applying Lemma (2.3), we get $\sum_{j}a_{j}(\Cal L_{t_{0}}.E_{j})/ 
\sum_{j}a_{j}m_{x_{0}}(E_{j}) \geq \min_{j}\{(\Cal L_{t_{0}}.E_{j})/
m_{x_{0}}(E_{j})\}$. Set $\min_{j}\{(\Cal L_{t_{0}}.E_{j})/
m_{x_{0}}(E_{j})\} = (\Cal L_{t_{0}}.E_{1})/m_{x_{0}}(E_{1})$. 
Now, combining all these together, we calculate 
$(\Cal L_{t_{0}}.E_{1})/m_{x_{0}}(E_{1}) \leq \sum_{j}a_{j}(\Cal L_{t_{0}}.
E_{j})/ 
\sum_{j}a_{j}m_{x_{0}}(E_{j}) \leq (\Cal L_{t_{0}}.C_{0})/m_{x_{0}}(C_{0}) 
\leq a$. 
Since $\epsilon(\Cal L_{t_{0}}, x_{0}) \leq 
(\Cal L_{t_{0}}.E_{1})/m_{x_{0}}(E_{1})$ (because $E_{1}$ 
is now irreducible and reduced), we obtain 
$\epsilon(\Cal L_{t_{0}}, x_{0}) \leq a$. Set $\Cal X(a)_{i} := 
\text{Im}(\text{pr}_{i, 1} : \Cal C(a)_{i} \rightarrow \Cal X)$, 
where $\text{pr}_{i, 1}$ is the natural evaluation morphism (from the first 
factor) given by $\Cal C(a)_{i} \ni (x, y) \mapsto x \in \Cal X$. (Remind 
that $\Cal C(a)_{i}$ is a subscheme of the universal family of 
$\Cal H^{0} \times_{\Cal B} \Cal H^{1}$.) Then, by $\epsilon(\Cal L_{t_{0}}, 
x_{0}) \leq a$ and by $[x_{0} \in C_{0}] \in \Cal H(a)_{i}$, 
we have $\Cal X(a)_{i} \subset \Cal X(a)$ and 
$\cup_{i = 1}^{I} \Cal X(a)_{i} \subset \Cal X(a)$. 
On the other hand, by the definition of $\Cal H(a)_{i}$, we have in apriori 
$\Cal X(a) \subset \cup_{i = 1}^{I} \Cal X(a)_{i}$. 
Therefore $\Cal X(a) = \cup_{i = 1}^{I} \Cal X(a)_{i}$. 
Since $f : \Cal X \rightarrow \Cal B$ 
and the natural morphisms $\Cal C(a)_{i} \rightarrow \Cal B$ are all proper, 
$\text{pr}_{i,1}$ are also proper. Hence, $\Cal X(a)_{i}$ are all Zariski 
closed in $\Cal X$, therefore, so is their finite union $\Cal X(a)$. \qed

\par
\vskip 20pt

Keiji Oguiso
\par
Math. Inst. d. Univ. Essen, D-45117 Essen, Germany; \par
Math. Sci. Univ. Tokyo, 153--8914 Tokyo, Japan \par
E-mail: mat9g0\@spi.power.uni-essen.de

\Refs
\widestnumber\key{CP}


\ref
\key B1 
\by  Th. Bauer
\paper Seshadri constants of quartic surfaces 
\jour  Math. Ann.
\vol  309
\yr  1997
\pages 475-481
\endref

\ref
\key B2 
\by  Th. Bauer
\paper Seshadri constants and periods of polarized abelian 
varieties, appendixed by Th. Bauer and T. Szemberg 
\jour  Math. Ann.
\vol  312
\yr  1998
\pages 607-623
\endref

\ref
\key B3 
\by  Th. Bauer
\paper Seshadri constants on algebraic surfaces 
\jour  Math. Ann.
\vol  313
\yr  1999
\pages 547-583
\endref

\ref
\key CP
\by F. Campana and T. Peternell
\paper Algebraicity of the ample cone of projective varieties 
\jour J. reine angew. Math.
\vol 407
\yr 1990
\pages 160--166
\endref

\ref
\key EL
\by L. Ein and R. Lazarsfeld
\paper Seshadri constants on smooth surfaces
\jour Ast\'erisque
\vol 218
\yr 1993
\pages 177-186
\endref

\ref
\key K
\by J. Koll\'ar
\paper Rational curves on algebraic varieties. 
A series of Mordern Surveys in Mathematics
\jour Springer-Verlag
\vol 32
\yr 1996
\endref

\ref
\key S
\by A. Steffens
\paper Remarks on Seshadri constants
\jour Math. Z.
\vol 227
\yr 1998
\pages 505-510
\endref

\endRefs
\enddocument